\documentclass[12pt]{article}
\newcommand{\di}{\displaystyle}

\usepackage{amsmath}
\usepackage{amssymb}
\usepackage{amsthm}
\newtheorem{theorem}{Theorem}
\newtheorem{lemma}{Lemma}
\newtheorem{definition}{Definition}
\newtheorem{corollary}{Corollary}

\begin{document}
\begin{center}
{\Large Some Rarita-Schwinger Type Operators}\\[7mm]

{Charles F. Dunkl, Junxia Li, John Ryan and Peter Van Lancker}
\end{center}

\begin{abstract}
 In this paper we study a generalization of the classical Rarita-Schwinger type operators and construct their fundamental solutions. We give some basic integral formulas related to these operators. We also establish that the projection operators appearing in the Rarita-Schwinger operators and the Rarita-Schwinger equations are conformally invariant. We further obtain the intertwining operators for other operators related to the Rarita-Schwinger operators under actions of the conformal group.
\end{abstract}


{\bf Classification:} Primary 30G35; Secondary 53C27


\section{Introduction}

\par In representation theory for $O(n)$ and $SO(n)$, one considers functions $f: U\longrightarrow  \mathcal{H}$$_{k}$ where U is a domain in
$\mathbb{R}^n$ and $\mathcal{H}$$_{k}$ is the space of harmonic polynomials homogeneous of degree $k$. Such spaces are invariant under actions of $O(n)$. If one refines to the covering group $Spin(n)$ of $SO(n)$, one replaces spaces of harmonic polynomials with spaces of homogeneous polynomial solutions to the Euclidean Dirac equation arising in Clifford analysis. See \cite{BDS}. Clifford analysis is the study of and applications of Dirac type operators. In this context the Rarita-Schwinger operators arise. See \cite{BSSV1, BSSV2, Va1, Va2, LRV1, LRV2, LR}. The Rarita-Schwinger operators are generalizations of the Dirac operator which in turn is a natural generalization of the Cauchy-Riemann operator. Rarita-Schwinger operators are also known as Stein-Weiss operators after \cite{SW}. We denote a Rarita-Schwinger operator by $R_k$, where $k=0, 1, \cdots, m, \cdots.$  When $k=0$ we have the Dirac operator.

\par  Here we start by constructing the Rarita-Schwinger operators and their fundamental solutions. This is based on the fundamental solution of the Dirac operator. Next, we give a summary of results on Rarita-Schwinger operators appearing in \cite{BSSV1}, giving detailed proofs and extending some of those results. We present a more detailed and alternative approach to that given in \cite{BSSV1}. This includes a Stokes' Theorem, Borel-Pompeiu Formula, Cauchy's Integral Formula and a Cauchy Transform. We also obtain intertwining operator for $R_k$ under actions of the conformal group, together with intertwining operators for the kernels to the Rarita-Schwinger operators, and the conformal invariance of Cauchy's Theorem and Cauchy's Integral Formula.
\par All of this ultimately helps to build the basics of Rarita-Schwinger type operators, including a theory of Rarita-Schwinger operators on examples of conformally flat manifolds. See for instance \cite{LRV1, LRV2, LR}.


\section{Preliminaries}

\par A Clifford algebra, $Cl_{n},$ can be generated from $\mathbb{R}^n$ by considering the
relationship $$\underline{x}^{2}=-\|\underline{x}\|^{2}$$ for each
$\underline{x}\in \mathbb{R}^n$.  We have $\mathbb{R}^n\subseteq Cl_{n}$. If $e_1,\ldots, e_n$ is an orthonormal basis for $\mathbb{R}^n$, then $\underline{x}^{2}=-\|\underline{x}\|^{2}$ tells us that $e_i e_j + e_j e_i= -2\delta_{ij}.$ Let $A=\{j_1, \cdots, j_r\}\subset \{1, 2, \cdots, n\}$ and $1\leq j_1< j_2 < \cdots < j_r \leq n$. An arbitrary element of the basis of the Clifford algebra can be written as $e_A=e_{j_1}\cdots e_{j_r}.$ Hence for any element $a\in Cl_{n}$, we have $a=\sum_Aa_Ae_A,$ where $a_A\in \mathbb{R}.$ For $a\in Cl_n$, we will need the following anti-involutions:
\begin{itemize}
\item Reversion:
 $$\tilde{a}=\sum_A(-1)^{|A|(|A|-1)/2}a_Ae_A,$$ where $|A|$ is the cardinality of $A$. In particular, $\widetilde{e_{j_1}\cdots e_{j_r}}=e_{j_r}\cdots e_{j_1}.$ Also $\widetilde{ab}=\tilde{b}\tilde{a}$ for $a, b \in Cl_n.$
\item Clifford conjugation:
 $$\bar{a}=\sum_A(-1)^{|A|(|A|+1)/2}a_Ae_A$$
satisfying $\overline{e_{j_1}\cdots e_{j_r}}=(-1)^r e_{j_r}\cdots e_{j_1}$ and $\overline{ab}= \bar{b}\bar{a}$ for $a, b \in Cl_n.$
\end{itemize}
For each $a=a_0+\cdots +a_{1\cdots n}e_1\cdots e_n\in Cl_n$ the scalar part of $\bar{a}a$ gives the square of the norm of $a,$ namely $a_0^2+\cdots +a_{1\cdots n}^2$\,.

The Pin and Spin groups play an important role in Clifford analysis.  The Pin group can be defined as
 $$Pin(n): =\{a\in Cl_n : a=y_1 \ldots y_p:
{y_1,\ldots , y_p}\in \mathbb{S}^{n-1}, p\in \mathbb{N}\}$$
and is clearly a group under multiplication in $Cl_n$.
\par Now suppose that $y\in \mathbb{S}^{n-1}\subseteq \mathbb{R}^n$. Look at $yxy=yx^{\parallel _y}y+yx^{\perp_y}y=-x^{\parallel _y}+x^{\perp_y}$ where $x^{\parallel _y}$ is the projection of $x$ onto $y$
and $x^{\perp_y}$ is perpendicular to $y$. So $yxy$ gives a reflection of $x$ in the $y$ direction. By the Cartan$-$Dieudonn\'{e} Theorem each $O \in O(n)$  is the composition of a finite number of reflections. If  $a=y_1\ldots y_p\in Pin(n)$, then $\tilde{a}:=y_p\ldots y_1$ and $ax\tilde{a}=O_a(x)$ for some $O_a\in O(n).$ Choosing $y_1, \ldots, y_p$ arbitrarily in $\mathbb{S}^{n-1}$,  we see that the group homomorphism $$\theta: Pin(n)\longrightarrow O(n): a\longmapsto O_a$$ with $a=y_1\ldots y_p$
and $O_a(x)=ax\tilde{a}$ is surjective. Further $-ax(-\tilde{a})=ax\tilde{a}$, so $1, -1\in ker(\theta)$. In fact $ker(\theta)=\{\pm 1\}.$
The Spin group is defined as
$$
Spin(n):=\{a\in Pin(n): a=y_1\ldots y_ p \mbox{ and } p \mbox{ even}\}
$$
and is a subgroup of $Pin(n)$. There is a group homomorphism
 $$\theta: Spin(n)\longrightarrow SO(n)$$ which is surjective with kernel $\{1, -1\}$. See \cite{P} for details.
\par The Dirac Operator in $\mathbb{R}^n$ is defined to be $$D :=\sum_{j=1}^{n} e_j \frac{\partial}{\partial x_j}.$$ Note $D^2=-\Delta_{n},$ where $\Delta_n$ is the Laplacian in $\mathbb{R}^n$.
\par Let $\mathcal{M}_k$ denote the space of $Cl_n-$ valued polynomials, homogeneous of degree $k$ and such that if $p_k\in$ $\mathcal{M}_k$ then $Dp_k=0.$ Such a polynomial is called a left monogenic polynomial homogeneous of degree $k$. Note if $h_k\in$ $\mathcal{H}_k,$ the space of $Cl_n-$ valued harmonic polynomials homogeneous of degree $k$, then $Dh_k\in$ $\mathcal{M}$$_{k-1}$. But $Dup_{k-1}(u)=(-n-2k+2)p_{k-1}(u),$
so $$\mathcal{H}_k=\mathcal{M}_k\bigoplus u\mathcal{M}_{k-1}, h_k=p_k+up_{k-1}.$$
This is the so-called Almansi-Fischer decomposition of $\mathcal{H}$$_k$. See \cite{BDS, R2}.
\par Note that if $Df(u)=0$ then $\bar{f}(u)\bar{D}=-\bar{f}(u)D=0$. So we can talk of right monogenic polynomials, homogeneous of degree $k$ and we obtain by conjugation a right Almansi-Fisher decomposition, $$\mathcal{H}_k=\overline{\mathcal{M}_k}\bigoplus \overline{\mathcal{M}}_{k-1}u,$$ where $\overline{\mathcal{M}_k}$ stands for the space of right monogenic polynomials homogeneous of degree $k$.


\section{\bf The Rarita-Schwinger Operator $R_k$}
Suppose $U$ is a domain in $\mathbb{R}^n$. Consider a function of two variables
$$f: U\times \mathbb{R}^n\longrightarrow Cl_n$$
such that for each $x\in U, f(x,u)$ is a left monogenic polynomial homogeneous of degree $k$ in $u$. Consider the action of the Dirac operator:
$$
D_xf(x,u)\,.
$$
As $Cl_n$ is not commutative then $D_xf(x,u)$ is no longer monogenic in $u$ but it is still harmonic and  homogeneous of degree $k$ in $u$.
So by the Almansi-Fischer decomposition, $D_xf(x,u)=f_{1,k}(x,u)+uf_{2,k-1}(x,u)$ where $f_{1,k}(x,u)$ is a left monogenic polynomial homogeneous of degree $k $ in $u$ and $f_{2,k-1}(x,u)$ is a left monogenic polynomial homogeneous of degree $k-1$ in $u$. Let $P_k$ be the left projection map
 $$P_k:  \mathcal{H}_k\rightarrow \mathcal{M}_k,$$ then
$R_kf(x,u)$ is defined to be $P_kD_xf(x,u)$. The left Rarita-Schwinger equation is defined to be
(see \cite{BSSV1})
$$
R_k f(x,u)=0.
$$
We also have a right projection $P_{k,r}: \mathcal{H}_k \rightarrow \overline{\mathcal{M}_k},$ and a right Rarita-Schwinger equation $f(x,u)D_xP_{k,r}=f(x,u)R_k=0.$ Since
$$
D_xf(x,u)=p_k(x,u)+up_{k-1}(x,u) \quad \mbox{and} \quad D_uup_{k-1}(x,u)=-(n+2k-2)p_{k-1}(x,u),
$$
we have $up_{k-1}(x,u)=-\di\frac{1}{n+2k-2}uD_uD_xf(x,u).$ Thus $(1-P_k)D_xf(x,u)=up_{k-1}(x,u)=-\di\frac{1}{n+2k-2}uD_uD_xf(x,u).$ Hence $$P_kD_xf(x,u)=\di\frac{1}{n+2k-2}uD_uD_xf(x,u)+D_xf(x,u)=(\di\frac{uD_u}{n+2k-2}+1)D_xf(x,u).$$ So we obtain that $P_k=(\di\frac{uD_u}{n+2k-2}+1)$ and $R_k=(\di\frac{uD_u}{n+2k-2}+1)D_x.$ See \cite{BSSV1}.

\par It is crucial to ask if there are any non-trivial solutions to this equation. First for any $k$-monogenic polynomial $p_k(u)$ we have trivially $R_kp_k(u)=0$. In particular the reproducing kernel of $\mathcal{M}_k$ is annihilated by $R_k$. We now produce a representation of this reproducing kernel. Consider the fundamental solution $G(u)=\di\frac{1}{\omega_n}\di\frac{-u}{\|u\|^n}$ to the Dirac operator $D$, where $\omega _n$ is the surface area of the unit sphere, $\mathbb{S}^{n-1}$.

Consider the Taylor series expansion of $G(v-u)$ and restrict to the $k$th order terms in $u_1,\ldots , u_n$ $(u=u_1e_1+\ldots +u_ne_n).$ These terms have as vector valued coefficients $$\di\frac{\partial ^k}{\partial {v_1}^{k_1}\ldots \partial {v_n}^{k_n}}G(v) \quad  (k_1+\ldots +k_n=k).$$
As $G(v)$ is a solution to the Dirac equation, $DG(v)=\sum_{i=1}^{n}e_j\di\frac{\partial G(v)}{\partial {v_j}}=0$, we can replace $\di\frac{\partial}{\partial{v_1}}$ by $-\sum_{j=2}^{n}e_{1}^{-1}e_j\di\frac{\partial}{\partial {v_j}}.$ Doing this each time $\di\frac{\partial}{\partial{v_1}}$ occurs and collecting like terms we obtain a finite series of polynomials homogeneous of degree $k $ in $u$ $$\di\sum_\sigma P_\sigma(u)V_\sigma(v)$$ where the summation is taken over all permutations of monogenic polynomials $(u_2-u_1e_1^{-1}e_2),\cdots, (u_n-u_1e_1^{-1}e_n)$, each term in the summation contains $j_2$ copies of $(u_2-u_1e_1^{-1}e_2),\cdots, j_n$ copies of $(u_n-u_1e_1^{-1}e_n),$ and
$$
P_\sigma(u)=\di\frac{1}{k!}\Sigma(u_{i_1}-u_1e_1^{-1}e_{i_1})\ldots(u_{i_k}-u_1e_1^{-1}e_{i_k}),
V_\sigma(v)=\di\frac{\partial^kG(v)}{\partial v_{2}^{j_2}\ldots \partial v_{n}^{j_n}}\,
$$
$j_2+\ldots +j_n=k,~~\mbox{and}~~ i_k\in\{2,\cdots,n\}.$ Here summation is taken over all permutations of the monomials without repetition. See \cite{BDS}. Note that this series is the sum of the $k-$th order terms in the Taylor expansion of $G(v-u)$ and consequently it is a vector.
 \par Now $\di\int_{\mathbb{S}^{n-1}}V_\sigma(u)uP_\mu(u)dS(u)=\delta_{\sigma,\mu}$ where $\delta_{\sigma,\mu}$ is the Kronecker delta and $\mu$ is a set of $n-1$ non-negative integers summing to $k$. See \cite{BDS}. Following \cite{BDS} it can be seen that the polynomial $P_\sigma$ is left monogenic and the set of all such polynomials, homogeneous of degree $k,$ forms a basis for the right $Cl_n$ module $\mathcal{M}_k.$
Consequently, the expression $$Z_k(u,v):=\di\sum_\sigma P_\sigma(u)V_\sigma(v)v$$ is the reproducing kernel of $\mathcal{M}_k$ with respect to integration over $\mathbb{S}^{n-1}$ (see \cite{BDS}).
Further as $Z_k(u,v)$ does not depend on $x$, $$R_kZ_k=0.$$
Note that $V_\sigma(v)v$ on $\mathbb{S}^{n-1}$ extends to $V_\sigma(-v^{-1})G(v)$ on $\mathbb{R}^n$ and this function is a right monogenic polynomial in $v$ and it is homogeneous of degree $k$. See \cite{BDS} and elsewhere.

We may ask if there are any solutions to $R_kf=0$ that depend on $x$. To do this we look at the interaction of the operator $R_k$ with conformal transformations.

\section{\bf Conformal transformations}
 We first establish the invariance properties of the projection operator $P_k$.
\subsection{\bf Conformal invariance of the projection $P_k$}
 Let $P_{k,w}$ and $P_{k,u}$ be the projections with respect to $w$ and $u$ respectively.
\subsubsection{\bf Orthogonal transformations}
\par Let $x=ay\tilde a,$ and $u=aw\tilde a.$
\begin{lemma}\quad $P_{k,w}\tilde af(ay\tilde a,aw\tilde a)=\tilde a P_{k,u}f(x,u),$ where $a\in Pin(n)$.\end{lemma}
\par  {\bf Proof} \quad  Let $f(x,u)=f_1(x,u)+uf_2(x,u),$ where $f_1(x,u)$ and $f_2(x,u)$ are monogenic polynomials homogeneous of degree $k$ and $k-1$ in $u.$  So $P_{k,u}f(x,u)=f_1(x,u)=f_1(ay\tilde a,aw\tilde a)$ and
$\tilde af(ay\tilde a,aw\tilde a)=\tilde af_1(ay\tilde a,aw\tilde a)+\tilde aaw\tilde a f_2(ay\tilde a,aw\tilde a)=\tilde af_1(ay\tilde a,aw\tilde a)\pm w\tilde a f_2(ay\tilde a,aw\tilde a).$
 \par Further as $\tilde af_1(ay\tilde a,aw\tilde a)$ and $\tilde af_2(ay\tilde a,aw\tilde a)$ are monogenic polynomials homogeneous of degree $k$ and $k-1$ in $w$ respectively, we have $P_{k,w}\tilde af(ay\tilde a,aw\tilde a)=\tilde af_1(ay\tilde a,aw\tilde a)=\tilde a P_{k,u}f(x,u).$ \qquad $\blacksquare$
\subsubsection{\bf Inversion} Let $x=y^{-1}, u=\di\frac{ywy}{\|y\|^2}.$
\begin{lemma} \quad $P_{k,w}\di\frac{y}{\|y\|^n}f(y^{-1},\di\frac{ywy}{\|y\|^2})=\di\frac{y}{\|y\|^n}P_{k,u}f(x,u).$\end{lemma}
\par {\bf Proof }\quad Since $f(x,u)=f_1(x,u)+uf_2(x,u),$ by substitution we have $$f(y^{-1},\di\frac{ywy}{\|y\|^2})=f_1(y^{-1},\di\frac{ywy}{\|y\|^2})+\di\frac{ywy}{\|y\|^2}f_2(y^{-1},\di\frac{ywy}{\|y\|^2}).$$
\par Now multiplying both sides of the above equation by $\di\frac{y}{\|y\|^n}$, one gets
$$\begin{array}{ll}\di\frac{y}{\|y\|^n}f(y^{-1},\di\frac{ywy}{\|y\|^2})=\di\frac{y}{\|y\|^n}f_1(y^{-1},\di\frac{ywy}{\|y\|^2})+\di\frac{y}{\|y\|^n}\di\frac{ywy}{\|y\|^2}f_2(y^{-1},\di\frac{ywy}{\|y\|^2})\\
\\
=\di\frac{y}{\|y\|^n}f_1(y^{-1},\di\frac{ywy}{\|y\|^2})-w\di\frac{y}{\|y\|^n}f_2(y^{-1},\di\frac{ywy}{\|y\|^2}).
\end{array}$$
\par Now Let $P_{k,w}$ act on the previous equation. We have $$P_{k,w}\di\frac{y}{\|y\|^n}f(y^{-1},\di\frac{ywy}{\|y\|^2})=\di\frac{y}{\|y\|^n}f_1(y^{-1},\di\frac{ywy}{\|y\|^2})=\di\frac{y}{\|y\|^n}f_1(x,u)=\di\frac{y}{\|y\|^n}P_{k,u}f(x,u),$$
which follows from the facts that $\di\frac{y}{\|y\|^n}f_1(y^{-1},\di\frac{ywy}{\|y\|^2})$ and $\di\frac{y}{\|y\|^n}f_2(y^{-1},\di\frac{ywy}{\|y\|^2})$ are monogenic and homogeneous of degree $k$ and $k-1$ in $w.$ \qquad $\blacksquare$
\subsubsection{\bf Translations} Let $x=y+a, a\in \mathbb{R}^n.$ In order to keep the homogeneity of $f(x,u)$ in $u$, $u$ does not change under translation. So we have
\begin{lemma} \quad $P_kf(x,u)=P_kf(y+a,u),$ where $x=y+a.$\end{lemma}
\subsubsection{\bf  Dilations}\quad Let $x=\lambda y,$ where $\lambda \in \mathbb{R}^{+}.$ It is obvious to observe that $P_k$ is invariant under dilation.
\begin{lemma} \quad $P_kf(x,u)=P_kf(\lambda y,u),$ where $x=\lambda y.$\end{lemma}
\par Hence $P_k$ is conformally invariant.\\

\par Ahlfors \cite{A} and Vahlen \cite{V} show that given a M\"{o}bius transformation $y=\phi(x)$ on $\mathbb{R}^n\bigcup\{\infty\}$ it can be expressed as $y=(ax+b)(cx+d)^{-1}$ where $a, b,c,d\in Cl_n$ and satisfy the following conditions:
\begin{enumerate}
 \item $a, b, c, d$ are all products of vectors in $\mathbb{R}^n.$
 \item $a\tilde b, c\tilde d, \tilde bc, \tilde da \in\mathbb{R}^n.$
 \item $a\tilde d-b\tilde c=\pm1.$
\end{enumerate}
When $c=0, \phi(x)=(ax+b)(cx+d)^{-1}=axd^{-1}+bd^{-1}=\pm ax\tilde a+bd^{-1}.$ Now assume $c\neq 0,$ then $\phi(x)=(ax+b)(cx+d)^{-1}=ac^{-1}\pm(cx\tilde{c}+d\tilde{c})^{-1},$ this is the so-called  Iwasawa decomposition. Using this notation and the conformal weights, $f(\phi(x))$ is changed to $J(\phi,x)f(\phi(x)),$ where $J(\phi,x)=\di\frac{\widetilde{cx+d}}{\|cx+d\|^n}$. Note when $\phi(x)=x+a$ then $J(\phi,x)\equiv 1.$ Now using the Iwasawa decomposition, we get the following result:\\
\begin{theorem} \quad $P_{k,w}J(\phi,x)f(\phi(x),\di\frac{\widetilde{(cx+d)}w(cx+d)}{\|cx+d\|^2})=J(\phi,x)P_{k,u}f(\phi(x),u),$ where $u=\di\frac{\widetilde{(cx+d)}w(cx+d)}{\|cx+d\|^2}$ and where $P_{k,w}$ and $P_{k,u}$ are the projections with respect to $w$ and $u$ respectively.\end{theorem}

Note that if the M\"{o}bius transformation is either translation or dilation then $u=w.$ This explains why in Lemma 4 the term $u$ is not multiplied by $\lambda.$
\par Lemmas 1 and 2 and Theorem 1 establish intertwining relationships for the projection operator, $P_k,$ under actions of the conformal group.


\subsection{\bf Conformal invariance of the Rarita-Schwinger operator $R_k$}

\par Now let us establish the intertwining operators for $R_k$ and the conformal invariance of the equation $R_kf=0.$ Let $R_{k,u}$ and $R_{k,w}$ be the Rarita-Schwinger operators with respect to $u$ and $w$ respectively.

We will need the following. If we have the M\"{o}bius transformation $y=\phi(x)$ and $D_x$ is the Dirac operator with respect to $x$ and $D_y$ is the Dirac operator with respect to $y$ then $D_x=J_{-1}(\phi,x)^{-1}D_yJ_1(\phi,x)$, where $J_{-1}(\phi,x)=\di\frac{cx+d}{\|cx+d\|^{n+2}}$ and $J_1(\phi,x)=J(\phi,x)=\di\frac{\widetilde{cx+d}}{\|cx+d\|^n}.$ See \cite{R1}.
 \subsubsection{\bf Orthogonal transformations $O\in O(n), a\in Pin(n)$}
\begin{theorem} \quad If $x=ay\tilde{a}, u=aw\tilde{a}$, then $aR_{k,u}f(x,u)=R_{k,w} \tilde{a}f(ay\tilde{a}, aw\tilde{a}).$
\end{theorem}
{\bf Proof} \quad $$R_{k,u}f(x,u)=P_{k,u}D_xf(x,u)=P_{k,u}a^{-1}D_y\tilde af(ay\tilde a,u)$$
Therefore, by Lemma 1 $$aP_{k,u}a^{-1}D_y\tilde af(ay\tilde a,u)=P_{k,w}aa^{-1}D_y\tilde af(ay\tilde a,aw\tilde{a})=R_{k,w} \tilde{a}f(ay\tilde{a}, aw\tilde{a}).\quad \blacksquare$$
\par In fact, Theorem 2 tells us that if $R_kf(x,u)=0$ then $R_k \tilde{a}f(ay\tilde{a}, aw\tilde{a})=0.$

\subsubsection{Inversion}
Let $x=y^{-1},(=\di\frac{-y}{\|y\|^2})$.
\begin{theorem}\quad Set $u=\di\frac{ywy}{\|y\|^2}$, then $\di\frac{y}{\|y\|^{n+2}}R_{k,u}f(x,u)=R_{k,w}G(y)f(y^{-1},\di\frac{ywy}{\|y\|^2}).$ \end{theorem}
{\bf Proof} \quad $$R_{k,u}f(x,u)=P_{k,u}D_xf(x,u)=P_{k,u}G_{-1}(y)^{-1}D_yG(y)f(y^{-1},u),$$ where $G_{-1}(y)=y\|y\|^n.$

Therefore by Lemma 2, $$\begin{array}{ll}G_{-1}(y)P_{k,u}G_{-1}(y)^{-1}D_yG(y)f(y^{-1},u)\\
\\
=P_{k,w}G_{-1}(y)G_{-1}(y)^{-1}D_yG(y)f(y^{-1},\di\frac{ywy}{\|y\|^2})=R_{k,w}G(y)f(y^{-1},\di\frac{ywy}{\|y\|^2}).\quad \blacksquare\end{array}$$

Consequently, if $R_kf(x,u)=0$, then $R_kG(y)f(y^{-1},\di\frac{ywy}{\|y\|^2})=0.$
\subsubsection{\bf Dilations} Let $x=\lambda y, \lambda \in \mathbb{R}^+.$ $R_kf(x,u)=R_kf(\lambda y,u)$ and if $R_kf(x,u)=0$ then $R_kf(\lambda y,u)=0.$
\subsubsection{\bf Translations} Let $x=y+ a , \quad a\in \mathbb{R}^n.$ In order to preserve homogeneity of polynomials in $u$, $f(x,u)$ is transformed under a translation by $a$ to $f(y+a, u)$ (Note: otherwise the action of the Vahlen matrices is not correct).   So $R_kf(x,u)=R_kf(y+a,u)$ and $R_kf(x,u)=0$ implies $R_kf(y+a, u)=0,$ where $x=y+a.$

Now using the Iwasawa decomposition of $(ax+b)(cx+d)^{-1}$, we obtain intertwining operators for $R_k:$
\begin{theorem}$$R_{k,x,w}J_1(\phi,x)\psi(\phi(x),\di\frac{\widetilde{(cx+d)}w(cx+d)}{\|cx+d\|^2})=J_{-1}(\phi,x)R_{k,y,u}\psi(y,u),$$ where $y=\phi(x), u=\di\frac{\widetilde{(cx+d)}w(cx+d)}{\|cx+d\|^2}, R_{k,x,w}=P_{k,w}D_x$ and $R_{k,y,u}=P_{k,u}D_y.$\end{theorem}
\par Consequently, we obtain that $R_kf(x,u)=0$ implies $$R_kJ(\phi,x)f(\phi(x),\di\frac{\widetilde{(cx+d)}w(cx+d)}{\|cx+d\|^2})=0,$$ where $u=\di\frac{\widetilde{(cx+d)}w(cx+d)}{\|cx+d\|^2}.$ For this last formula see also \cite{BSSV1}.


\section{\bf A  Kernel for $R_k$ and Some Basic Integral Formulas}
\par  Now applying inversion from the left we obtain that if $R_kZ_k(u,v)=0$ then \\$R_kG(x)Z_k(\di\frac{xux}{\|x\|^2},v)=0$. That is, $$F_k(x,u,v)= \di\frac{x}{\|x\|^n}Z_k(\frac{xux}{\|x\|^2}, v)=
\di\frac{x}{\|x\|^{n+2k}}Z_k(xux, v)$$ is a non-trivial solution to $R_kf(x, u)=0$ on $\mathbb{R}^n\setminus{\{0\}}.$ Note that this function is monogenic in $u$.

Similarly, applying inversion from the right we obtain that $$ Z_k(u,\frac{xvx}{\|x\|^2} )\di\frac{x}{\|x\|^n}= Z_k(u,xvx)\di\frac{x}{\|x\|^{n+2k}}$$ is a non-trivial solution to $f(x, v)R_k=0$ on $\mathbb{R}^n\setminus{\{0\}}.$
In fact, \cite{BSSV1}, this function is $F_k(x,u,v),$ and, up to a constant, $F_k(x,u,v)$ is the fundamental solution of $R_k$. One proof of this statement is given in the following. Let $\mathcal{M}_k$ denote the space of left monogenic polynomials homogeneous of degree $k$ and suppose $a\in Pin(n),$ then $$l_{\tilde a}:\mathcal{M}_k\rightarrow \mathcal{M}_k\quad : f(u)\to \tilde{a}f(au\tilde a)$$ is an isomorphism. Similarly, let $\overline{\mathcal{M}_k}$ denote the space of right monogenic polynomials homogeneous of degree $k$ then
$$r_{\tilde a}: \overline{\mathcal{M}_k}\rightarrow\overline{\mathcal{M}_k}\quad : f(u)\to f(au\tilde a)a$$ is also an isomorphism. Using these isomorphisms it may be seen that for each $a\in Pin(n)$ then $\pm \tilde aZ_k(au\tilde a,av\tilde a)a$ is also the reproducing kernel for $\mathcal{M}_k$. We choose the plus sign when $a\in Spin(n)$ and the minus sign when $a\in Pin(n)\setminus Spin(n).$ Now consider a non-zero vector $x\in \mathbb{R}^n$, then $\di\frac{x}{\|x\|}\in Pin(n)$. So we have
$$Z_k(u,v)=-\di\frac{x}{\|x\|}Z_k(\di\frac{xux}{\|x\|^2},\di\frac{xvx}{\|x\|^2})\di\frac{x}{\|x\|},$$ that is,
$-\di\frac{x}{\|x\|}Z_k(\di\frac{xux}{\|x\|^2},\di\frac{xvx}{\|x\|^2})\di\frac{x}{\|x\|}$ is also the reproducing kernel for $\mathcal{M}_k$.
Now we look at the fundamental solution of $R_k$ which has the representation $Z_k(u,\di\frac{xvx}{\|x\|^2})\di\frac{x}{\|x\|^n}$. Then using the previous equality, we get $$\begin{array}{ll}
Z_k(u,\di\frac{xvx}{\|x\|^2})\di\frac{x}{\|x\|^n}
=-\di\frac{x}{\|x\|}Z_k(\di\frac{xux}{\|x\|^2},v)\di\frac{x}{\|x\|}\di\frac{x}{\|x\|^n}=
\di\frac{x}{\|x\|^n}Z_k(\di\frac{xux}{\|x\|^2},v).\end{array}$$

\par Further suppose $\mu$ is a $Cl_n$ valued measure on $\mathbb{R}^n$ with compact support, $[\mu].$ It follows for suitable choices of $\mu$ the integral $\di\int_{[\mu]}F_k(x,u,v)d\mu$ defines a solution to $R_kf=0$ on $(\mathbb{R}^n\setminus[\mu])\times \mathbb{R}^n.$

 \subsection{\bf Stokes' Theorem}
We first build Stokes' Theorem for the Rarita-Schwinger operator. This is based on Stokes' Theorem for the Dirac operator.
 \begin{theorem} $($Stokes' Theorem for the Dirac operator, \cite{BDS}$)$ \\
 Let $\Omega$ and $\Omega'$ be domains in $\mathbb{R}^n$ and suppose the closure of $\Omega$ lies in $\Omega'$. Further suppose the closure of $\Omega$ is compact and $\partial\Omega$ is piecewise smooth. Let $f, g \in C^1(\Omega', Cl_n).$ Then $$\int_{\partial\Omega}g(x,u)d\sigma_xf(x,u)=\int_\Omega[\left(g(x,u)D_x\right)f(x,u)+g(x,u)\left(D_xf(x,u)\right)]dx^n,$$
where $dx^n = dx_1\wedge\cdots \wedge dx_n$, $d\sigma_x=n(x)d\sigma(x),$ $\sigma$ is scalar Lebesgue measure on $\partial\Omega$ and $n(x)$ is unit outer normal vector to $\partial\Omega.$ We may write $n(x)$ as $\di\sum_{i=1}^{n}n_i(x)e_i,$ where $n_i(x)$ are scalar-valued functions. $g(x,u)D_x$ means $D_x$ acts from the right on $g(u,x).$
 \end{theorem}

\begin{definition}  \quad For any $Cl_n-$valued polynomials $P(u), Q(u)$, the inner product $(P(u), Q(u))_u$ with respect to $u$ is given by $$(P(u), Q(u))_u=\di\int_{\mathbb{S}^{n-1}}P(u)Q(u)dS(u).$$
\end{definition}
This inner product differs slightly from the Fischer inner product in \cite{BSSV1}. There the inner product is $\di\int_{\mathbb{S}^{n-1}}\overline{R}(u)Q(u)dS(u)$ for a $Cl_n$ valued polynomial $R(u).$ If we place $R(u)=\overline{P}(u)$ we see that, as the conjugation$^-$ is an isomorphism, the two inner products are equivalent. For any $p_k \in \mathcal{M}_k,$ one obtains (see \cite{BDS})
$$
p_k(u)=(Z_k(u,v), p_k(v))_v=\int_{\mathbb{S}^{n-1}}Z_k(u,v)p_k(v)dS(v).
$$
Using Stokes' Theorem for the Dirac operator, we can obtain the basic formulas related to the Rarita-Schwinger operators.\\
\begin{lemma}
Suppose $p_k$ is a left monogenic polynomial homogeneous of degree $k$ and $p_{k-1}$ is a left monogenic polynomial homogeneous of degree $k-1$ then
$$\di\int_{\mathbb{S}^{n-1}}\tilde{p}_{k-1}(u)up_k(u)dS(u)=0.$$
\end{lemma}
Outline Proof: As we are integrating over the unit sphere the previous integral can be written as
$$\di\int_{\mathbb{S}^{n-1}}\tilde{p}_{k-1}(u)n(u)p_k(u)dS(u).$$
By the Clifford-Cauchy Theorem \cite{BDS} this integral vanishes. $\blacksquare$

We now have
\begin{theorem} (Rarita-Schwinger Stokes' Theorem)\cite{BSSV1} Let $\Omega'$ and $\Omega$ be as in Theorem 5. Then for $f, g \in C^1(\Omega',$$\cal{M}$$_k)$, we have
$$\begin{array}{ll}
\di\int_{\partial\Omega}\left(g(x,u)d\sigma_xf(x,u)\right)_u\\
\\
=\di\int_\Omega(g(x,u)R_k, f(x,u))_udx^n+\di\int_\Omega(g(x,u), R_kf(x,u))_udx^n.
\end{array}$$
Further
$$\begin{array}{ll}
\di\int_{\partial\Omega}\left(g(x,u)d\sigma_xf(x,u)\right)_u\\
\\
=\di\int_{\partial\Omega}\left(g(x,u), P_kd\sigma_xf(x,u)\right)_u\\
\\
=\di\int_{\partial\Omega}\left(g(x,u)d\sigma_xP_k, f(x,u)\right)_u.
\end{array}$$\end{theorem}
Outline Proof: The first identity is obtained by first applying Stokes' Theorem to the integral
$\di\int_{\partial\Omega}\left(g(x,u)d\sigma_xf(x,u)\right)_u$ to obtain
$$\di\int_\Omega(g(x,u)D, f(x,u))_u+(g(x,u), Df(x,u))_udx^n.$$

Both $g(x,u)D$ and $Df(x,u)$ have an Almansi-Fischer decomposition with respect to $u$. So applying Lemma 5 with respect to $u$ and Definition 1 and these Almansi-Fischer decompositions give the result.

The second collection of identities again arise by applying the Almansi-Fischer decomposition $d\sigma_xf(x,u)$ and $g(x,u)d\sigma_x$ respectively, and then applying Definition 1 and Lemma 5 with respect to $u$. $\blacksquare$

Now if both $f(x,u)$ and $g(x,u)$ are solutions of $R_k$, then we have the following result:

\begin{corollary} $($Cauchy's Theorem$)$ \\
If $R_kf(x,u)=0$ and $g(x,u)R_k=0$ for $f, g \in C^1(\Omega',$$\cal{M}$$_k)$, then $$\di\int_{\partial\Omega}\left(g(x,u), P_kd\sigma_xf(x,u)\right)_u=0.$$\end{corollary}
Let $S$ be a hypersurface in $\mathbb{R}^n$ and $y=\phi(x)=(ax+b)(cx+d)^{-1}.$ Now look at Cauchy's Theorem:
\begin{eqnarray*}
0&=&\di\int_{S}\left(g(y,u), P_kd\sigma_yf(y,u)\right)_u=\di\int_{S}\left(g(y,u), P_kn(y)f(y,u)\right)_ud\sigma(y)\\
&=&\di\int_{\phi^{-1}(S)}\left(g(\phi(x),u), P_k\tilde J(\phi,x)n(x)J(\phi,x)f(\phi(x),u)\right)_ud\sigma(x)\\
&=&\di\int_{\phi^{-1}(S)}\int_{\mathbb{S}^{n-1}}g(\phi(x),u)P_{k,u}\tilde J(\phi,x)n(x)J(\phi,x)f(\phi(x),u)dS(u)d\sigma(x).
\end{eqnarray*}
Set $u=\di\frac{\widetilde{(cx+d)}w(cx+d)}{\|cx+d\|^2},$ since $P_{k,u}$ can interchange with $\tilde J(\phi,x)$, the previous equation equals
$$\begin{array}{ll}0=\di\int_{\phi^{-1}(S)}\int_{\mathbb{S}^{n-1}}g(\phi(x),\di\frac{\widetilde{(cx+d)}w(cx+d)}{\|cx+d\|^2})\tilde J(\phi,x)P_{k,w}n(x)J(\phi,x)f(\phi(x),\\
\\
\di\frac{\widetilde{(cx+d)}w(cx+d)}{\|cx+d\|^2}) dS(w)d\sigma(x)\\
\\
=\di\int_{\phi^{-1}(S)}(g(\phi(x),\di\frac{\widetilde{(cx+d)}w(cx+d)}{\|cx+d\|^2})\tilde J(\phi,x),P_{k,w}d\sigma_xJ(\phi,x)\\
\\
f(\phi(x),\di\frac{\widetilde{(cx+d)}w(cx+d)}{\|cx+d\|^2}))_w.
\end{array}$$

Therefore, Cauchy's Theorem is conformally invariant under M\"{o}bius transformations.

We now wish to introduce the Borel-Pompeiu Theorem from \cite{BSSV1}. First we will need:
\begin{lemma}
Suppose $h_k: \mathbb{R}^n\to Cl_n$ is a harmonic polynomial homogeneous of degree $k$ and $n>2$. Suppose $u\in\mathbb{S}^{n-1}$ then
$$\di\frac{1}{\omega_n}\di\int_{\mathbb{S}^{n-1}}h_k(xux)dS(x)=c_kh_k(u),$$ where $c_k=\di\frac{n-2}{n-2+2k}$.\end{lemma}

The proof follows from \cite{DX} (Proposition 2.2.3 on page 34).

{\bf Proof}\quad Fix $u\in \mathbb{S}^{n-1},$ suppose $h_k: \mathbb{R}^n\to Cl_n$ is a harmonic polynomial homogeneous of degree $k$. Now compute
$$I(h_k)=\di\frac{2}{\omega_n}\di\int_{\left\langle x,u\right\rangle >0}h_k(u-2\left\langle x,u\right\rangle x)dS(x),$$ note $xux=u-2\left\langle x,u\right\rangle x$ is the reflection of $u$ in the mirror $x^{\perp}$ (where $\di\int_{\left\langle x,u\right\rangle >0}dS(x)=\di\frac{1}{2}\int_{\mathbb{S}^{n-1}}dS(x)=\di\frac{\omega_n}{2}$, and $\|x\|=1$). Since $xux=u-2\left\langle x,u\right\rangle x$ is invariant under $x\to -x$ we have $$I(h_k)=\di\frac{1}{\omega_n}\di\int_{\mathbb{S}^{n-1}}h_k(u-2\left\langle x,u\right\rangle x)dS(x).$$
By rotation invariance assume $u=(0,\cdots,0,1).$ Any harmonic homogeneous polynomial of degree $k$ has an expression
$$h_k(y)=\di\sum_{j=0}^{k}\|y\|^{k-j}P_{k-j}^{j+n/2-1}\left(\di\frac{y_n}{\|y\|}\right)h_{k_j}(y_1,\cdots,y_{n-1}),$$
where $h_{k_j}$ is harmonic and homogeneous of degree $j$ and the normalized Gegenbauer polynomial is
$$P_m^{\lambda}(t)=\di\sum_{i=0}^{m}\di\frac{(-m)_i(m+2\lambda)_i}{(\lambda+\frac{1}{2})_ii!}\left(\di\frac{1-t}{2}\right)^i.$$
In the coordinate system $y=(y'\sin\theta  ,\cos\theta)$ with $0\leq\theta\leq\pi$ and $y'\in \mathbb{S}^{n-2}$ the integral is
$$\di\int_{\mathbb{S}^{n-1}}h_k(y)dS(x)=c_n'\di\int_0^{\pi}\sin^{n-2}\theta d\theta \di\int_{\mathbb{S}^{n-2}}h_k(y'\sin\theta  ,\cos\theta)dS_{n-2}(y').$$
Set $u=(0,\cdots,0,1)$ and $x=(x'\sin\theta  ,\cos\theta)$ with $x'\in \mathbb{S}^{n-2},$ then
$$u-2\left\langle x,u\right\rangle x=(-2x'\cos\theta\sin\theta  ,1-2\cos^2\theta).$$
Thus
$$\begin{array}{ll}
\di\int_{\mathbb{S}^{n-1}}h_k(u-2\left\langle x,u\right\rangle x)dS(x)\\
\\
=\di\sum_{j=0}^{k}c_n'\di\int_0^{\pi}P_{k-j}^{j+n/2-1}(1-2\cos^2\theta)\sin^{n-2}\theta d\theta \di\int_{\mathbb{S}^{n-2}}h_{k_j}(-2x'\cos\theta \sin\theta  )dS_{n-2}(x')\\
\\
=\di\sum_{j=0}^{k}c_n'\di\int_0^{\pi}P_{k-j}^{j+n/2-1}(1-2\cos^2\theta)(-2\cos\theta \sin\theta)^j\sin^{n-2}\theta d\theta \di\int_{\mathbb{S}^{n-2}}h_{k_j}( x')dS_{n-2}(x')\\
\\
=h_{k_0}c_n'\di\int_0^{\pi}P_{k}^{n/2-1}(1-2\cos^2\theta)\sin^{n-2}\theta d\theta.
\end{array}$$
The integrals equal zero for $j>0,$ by the orthogonality property of harmonics. Next set $t=\cos\theta$ and $dt=-\sin\theta d\theta,$ then
\begin{eqnarray*}
\di\int_0^{\pi}P_{k}^{\lambda}(1-2\cos^2\theta)\sin^{n-2}\theta d\theta &=& \di\int_{-1}^1P_k^{\lambda}(1-2t^2)(1-t^2)^{\lambda-\frac{1}{2}}dt\\
&=& \di\sum_{i=0}^{k}\di\frac{(-k)_i(k+2\lambda)_i}{(\lambda+\frac{1}{2})_ii!}\di\int_{-1}^1 t^{2i}(1-t^2)^{\lambda-\frac{1}{2}}dt\\
&=& \di\sum_{i=0}^{k}\di\frac{(-k)_i(k+2\lambda)_i}{(\lambda+\frac{1}{2})_ii!}B(i+\frac{1}{2},\lambda+\frac{1}{2})\\
&=& \di\frac{\Gamma(\frac{1}{2})\Gamma(\lambda+\frac{1}{2})}{\Gamma(\lambda+1)}\di\sum_{i=0}^{k}\di\frac{(-k)_i(k+2\lambda)_i(\frac{1}{2})_i}{(\lambda+\frac{1}{2})_ii!(\lambda+1)_i}.
\end{eqnarray*}
The Saalsch\"{u}tz summation formula (for $-k+a+b+1=c+d$) is
$$
{}_3 F_2\left(\begin{matrix}-k, a, b\\c, d\end{matrix};1\right)=\di\frac{(c-a)_k(d-a)_k}{(c)_k(d)_k}.$$
\begin{eqnarray*}\di\int_0^{\pi}P_{k}^{\lambda}(1-2\cos^2\theta)\sin^{n-2}\theta d\theta &
=&\di\frac{\Gamma(\frac{1}{2})\Gamma(\lambda+\frac{1}{2})}{\Gamma(\lambda+1)}\di\frac{(\lambda)_k(\lambda+\frac{1}{2})_k}{(\lambda+\frac{1}{2})_k(\lambda+1)_k}\\
&=&\di\frac{\Gamma(\frac{1}{2})\Gamma(\lambda+\frac{1}{2})}{\Gamma(\lambda+1)}\di\frac{\lambda}{\lambda+k}.
\end{eqnarray*}
The normalizing constant is (now $\lambda=\frac{n}{2}$)$c_n'=\di\frac{\Gamma(\frac{n}{2})}{\Gamma(\frac{1}{2})\Gamma(\frac{n-1}{2})}$ and
$$\di\frac{1}{\omega_n}\di\int_{\mathbb{S}^{n-1}}h_k(u-2\left\langle x,u\right\rangle x)dS(x)=\di\frac{\frac{n}{2}-1}{\frac{n}{2}-1+k}h_{k_0}=c_kh_{k_0},$$ and $h_k(0,0,\cdots,1)=h_{k_0},$ where $c_k=\di\frac{n-2}{n-2+2k}$.\quad $\blacksquare$

We now introduce $E_k(x,u,v):=\di\frac{1}{\omega_n c_k}F_k(x,u,v).$
\begin{theorem} \cite{BSSV1} (Borel-Pompeiu Theorem) Let $\Omega'$ and $\Omega$ be as in Theorem 5 and $y\in \Omega.$ Then for $f \in C^1(\Omega',\mathcal{M}_k)$
$$\begin{array}{ll}
f(y,u)=\di\int_{\partial\Omega}\left(E_k(x-y,u,v), P_{k,v}d\sigma_xf(x,v)\right)_v-\di\int_\Omega(E_k(x-y,u,v),R_{k,v}f(x,v))_vdx^n.
\end{array}$$ \end{theorem}
\par{\bf Proof}\qquad Here we will use the representation $$E_k(x-y,u,v)=\di\frac{1}{\omega_n c_k}Z_k(u,\di\frac{(x-y)v(x-y)}{\|x-y\|^2})\di\frac{x-y}{\|x-y\|^n},$$ and $R_{k,v}$ is the Rarita-Schwinger operator with respect to $v.$
 Consider a ball $B(y,r)$ centered at $y$ with radius $r$ such that $\overline{B(y,r)}\subset \Omega$. By Stokes' Theorem, we have
$$\begin{array}{ll}
\di\int_\Omega(E_k(x-y,u,v),R_{k,v}f(x,v))_vdx^n\\
\\
=\di\int_{\Omega\setminus{B(y,r)}}(E_k(x-y,u,v),R_{k,v}f(x,v))_vdx^n+\di\int_{B(y,r)}(E_k(x-y,u,v),R_{k,v}f(x,v))_vdx^n.
\end{array}$$
The second integral tends to zero as $r$ tends to zero. This follows from the degree of homogeneity of $E_k(x-y,u,v)$. Now applying Stokes' Theorem to the first integral, one gets
$$\begin{array}{ll}\di\int_{\Omega\setminus{B(y,r)}}(E_k(x-y,u,v),R_{k,v}f(x,v))_vdx^n\\
\\
=\di\int_{\partial\Omega}(E_k(x-y,u,v),P_{k,v}d\sigma_xf(x,v))_v-\di\int_{\partial {B(y,r)}}(E_k(x-y,u,v),P_{k,v}d\sigma_xf(x,v))_v.
 \end{array}$$
Now let us look at $$\begin{array}{ll}\di\int_{\partial {B(y,r)}}(E_k(x-y,u,v),P_{k,v}d\sigma_xf(x,v))_vdx^n=
\di\int_{\partial {B(y,r)}}(E_k(x-y,u,v),P_{k,v}d\sigma_xf(y,v))_v\\
\\
\qquad \qquad \qquad \qquad \qquad \qquad+\di\int_{\partial {B(y,r)}}(E_k(x-y,u,v),P_{k,v}d\sigma_x[f(x,v)-f(y,v)])_v.
 \end{array}$$
 Since the second integral on the right hand side tends to zero as $r$ goes to zero, we only need to deal with the first integral.
 $$\begin{array}{llll}\di\int_{\partial {B(y,r)}}(E_k(x-y,u,v),P_{k,v}d\sigma_xf(y,v))_v\\
 \\
 =\di\int_{\partial{B(y,r)}}\int_{\mathbb{S}^{n-1}}E_k(x-y,u,v)P_{k,v}d\sigma_xf(y,v)dS(v)\\
 \\
 =\di\int_{\partial{B(y,r)}}\int_{\mathbb{S}^{n-1}}\di\frac{1}{\omega_n c_k}Z_k\left(u,\di\frac{(x-y)v(x-y)}{\|x-y\|^2}\right)\di\frac{x-y}{\|x-y\|^n}P_{k,v}n(x)f(y,v)dS(v)d\sigma(x),
 \end{array}$$ where $n(x)$ is the unit outer normal vector and $d\sigma(x)$ is the scalar measure on $\partial B(y,r).$ Now $n(x)$ here is $\di\frac{y-x}{\|x-y\|}.$  Hence the previous integral becomes
 $$\begin{array}{ll}
 \di\frac{1}{\omega_nc_k}\di\int_{\partial{B(y,r)}}\int_{\mathbb{S}^{n-1}}Z_k\left(u,\di\frac{(x-y)v(x-y)}{\|x-y\|^2}\right)\di\frac{x-y}{\|x-y\|^n}P_{k,v}\di\frac{y-x}{\|x-y\|}f(y,v)
 dS(v)d\sigma(x).\end{array}$$

 By equation $(1)$ this integral becomes
 $$\begin{array}{ll}\di\frac{1}{\omega_nc_k}\di\int_{\partial{B(y,r)}}\int_{\mathbb{S}^{n-1}}Z_k\left(u,\di\frac{(x-y)v(x-y)}{\|x-y\|^2}\right)\di\frac{x-y}{\|x-y\|^n}\di\frac{y-x}{\|x-y\|}f(y,v)dS(v)d\sigma(x)\\
 \\
=\di\frac{1}{\omega_n}c_k^{-1}\di\int_{\partial{B(y,r)}}\di\frac{1}{r^{n-1}}\int_{\mathbb{S}^{n-1}}Z_k(u,\di\frac{(x-y)v(x-y)}{\|x-y\|^2})f(y,v)d\sigma(x)dS(v)
\end{array}$$

By Lemma 6 this integral is equal to
$$\begin{array}{ll}
\di\int_{\mathbb{S}^{n-1}}Z_k(u,v)f(y,v)dS(v)=f(y,u).
\quad \blacksquare\end{array}$$

\begin{theorem} \cite{BSSV1}(Cauchy's Integral Formula) If $R_kf(x,v)=0,$ then for $y\in \Omega,$
$$\begin{array}{ll}
f(y,u)=\di\int_{\partial\Omega}\left(E_k(x-y,u,v), P_kd\sigma_xf(x,v)\right)_v\\

=\di\int_{\partial\Omega}\left(E_k(x-y,u,v)d\sigma_xP_{k,r} f(x,v)\right)_v.
\qquad \blacksquare \end{array}$$ \end{theorem}

\par We now show the conformal invariance of Cauchy's Integral Formula. We start with inversion.
\par Since $$\begin{array}{ll}x^{-1}-y^{-1}=-y^{-1}(x-y)x^{-1}=-x^{-1}(x-y)y^{-1}\\
\\
=\di\frac{-x}{\|x\|^2}(x-y)\di\frac{y}{\|y\|^2}=\di\frac{-y}{\|y\|^2}(x-y)\di\frac{x}{\|x\|^2},\end{array}$$

$$\begin{array}{ll}E_k(x^{-1}-y^{-1},u,v)=G(x^{-1}-y^{-1})Z_k\left(\di\frac{(x^{-1}-y^{-1})u(x^{-1}-y^{-1})}{\|x^{-1}-y^{-1}\|^2},v\right)\\
\\
=-G(y)^{-1}G(x-y)G(x)^{-1}Z_k\left(\di\frac{x(x-y)yuy(x-y)x}{\|x\|^2\|y\|^2\|x-y\|^2},v\right)\\
\\
=-G(y)^{-1}G(x-y)G(x)^{-1}\di\frac{-x}{\|x\|}Z_k\left(\di\frac{(x-y)yuy(x-y)}{\|y\|^2\|x-y\|^2},\di\frac{xvx}{\|x\|^2}\right)\di\frac{x}{\|x\|}\\
\\
=G(y)^{-1}G(x-y)Z_k\left(\di\frac{(x-y)u'(x-y)}{\|x-y\|^2},\di\frac{xvx}{\|x\|^2}\right)x\|x\|^{n-2}, ~~\mbox{set}~~ u'=\di\frac{yuy}{\|y\|^2}\\
\\
=-G(y)^{-1}G(x-y)Z_k\left(\di\frac{(x-y)u'(x-y)}{\|x-y\|^2},\di\frac{xvx}{\|x\|^2}\right)G(x)^{-1}\\
\\
=-G(y)^{-1}E_k(x-y,u',v')G(x)^{-1},
\end{array}$$ where $u'=\di\frac{yuy}{\|y\|^2}$ and $v'=\di\frac{xvx}{\|x\|^2}$.
\par Now consider 
$$\begin{array}{ll}E_k(ax\tilde a-ay\tilde a,u,v)=G(a(x-y)\tilde a)Z_k\left(\di\frac{a(x-y)\tilde aua(x-y)\tilde a}{\|a(x-y)\tilde a\|^2},v\right)\\
\\
=aG(x-y)\tilde aZ_k\left(\di\frac{a(x-y)\tilde aua(x-y)\tilde a}{\|x-y\|^2},v\right)\\
\\
=\pm aG(x-y)\tilde aaZ_k\left(\di\frac{\tilde aa(x-y)\tilde aua(x-y)\tilde aa}{\|x-y\|^2},\tilde ava\right)\tilde a
\end{array}$$
$$\begin{array}{ll}
=aG(x-y)Z_k\left(\di\frac{(x-y)u'(x-y)}{\|x-y\|^2},\tilde ava\right)\tilde a,~~\mbox{set}~~u'=\tilde aua\\
\\
=aE_k(x-y,u',v')\tilde a,
\end{array}$$ where $u'=\tilde aua$ and $v'=\tilde ava$.
\par Using the Iwasawa decomposition, one gets $$E_k(\phi(x)-\phi(y),u,v)=J(\phi,y)^{-1}E_k(x-y,u',v')\tilde J(\phi,x)^{-1},$$ where $u'=\di\frac{\widetilde{(cy+d)}u(cy+d)}{\|cy+d\|^2}, v'=\di\frac{\widetilde{(cx+d)}v(cx+d)}{\|cx+d\|^2},$ and $\phi$ is the M\"{o}bius transformation.

\par Suppose $S$ is a smooth hypersurface lying in $\mathbb{R}^n.$ Let $x'=\phi(x)$ and $y'=\phi(y)$, now let us consider Cauchy's Integral Formula
$$\begin{array}{ll}
f(y',u)=\di\int_{S}(E_k(x'-y',u,v),P_kn(x')f(x',v))_vd\sigma(x')\\
\\
=\di\int_{S}\int_{\mathbb{S}^{n-1}}E_k(x'-y',u,v)P_kn(x')f(x',v)dS(v)d\sigma(x').
\end{array}$$
Thus, by the fact that $n(x')d\sigma(x')=\tilde J(\phi,x)n(x)J(\phi,x)d\sigma(x)$ we have
$$\begin{array}{ll}
f(\phi(y),u)=\di\int_{\phi^{-1}(S)}\int_{\mathbb{S}^{n-1}}J(\phi,y)^{-1}E_k(x-y,u',v')\tilde J(\phi,x)^{-1}P_{k,v}\tilde J(\phi,x)n(x)J(\phi,x)\\
\\
\qquad \qquad\qquad\qquad \qquad \qquad \qquad \qquad \qquad\qquad\qquad \qquad \qquad \qquad  f(\phi(x),v)dS(v)d\sigma(x)\\
\\
=\di\int_{\phi^{-1}(S)}\int_{\mathbb{S}^{n-1}}J(\phi,y)^{-1}E_k(x-y,u',v')\tilde J(\phi,x)^{-1}\tilde J(\phi,x)P_{k,v'}n(x)J(\phi,x)\\
\\
\qquad \qquad\qquad\qquad \qquad \qquad \qquad \qquad \qquad\qquad f(\phi(x),\di\frac{(cx+d)v'\widetilde{(cx+d)}}{\|cx+d\|^2})dS(v')d\sigma(x)
\end{array}$$
$$\begin{array}{ll}
=\di\int_{\phi^{-1}(S)}\int_{\mathbb{S}^{n-1}}J(\phi,y)^{-1}E_k(x-y,u',v')P_{k,v'}n(x)J(\phi,x)f(\phi(x),\di\frac{(cx+d)v'\widetilde{(cx+d)}}{\|cx+d\|^2})\\
\\
\qquad \qquad\qquad\qquad \qquad \qquad \qquad \qquad \qquad\qquad \qquad \qquad\qquad dS(v')d\sigma(x).
\end{array}$$
 Multiplying both sides of the previous equation by $J(\phi,y)$, we obtain
$$\begin{array}{ll}J(\phi,y)f(\phi(y),\di\frac{(cy+d)u'\widetilde{(cy+d)}}{\|cy+d\|^2})=\di\int_{\phi^{-1}(S)}\int_{\mathbb{S}^{n-1}}E_k(x-y,u',v')P_{k,v'}n(x)J(\phi,x)\\
\\
\qquad \qquad\qquad\qquad \qquad \qquad \qquad \qquad \qquad f(\phi(x),\di\frac{(cx+d)v'\widetilde{(cx+d)}}{\|cx+d\|^2})dS(v')d\sigma(x),\end{array}$$
where $u=\di\frac{(cy+d)u'\widetilde{(cy+d)}}{\|cy+d\|^2}$ and $v=\di\frac{(cx+d)v'\widetilde{(cx+d)}}{\|cx+d\|^2}.$
\par Therefore, Cauchy's Integral Formula is conformally invariant.\\

Now if the function $\psi$ has compact support in $\Omega$, then by the Borel-Pompeiu Theorem we have the following formula:
\begin{theorem} $\di\iint_{\mathbb{R}^n}-(E_k(x-y,u,v),R_k\psi(x,v))_vdx^n=\psi(y,u)$ for each $\psi\in C_0^{\infty}(\mathbb{R}^n,\mathcal{M}_k)$.
\end{theorem}

Similarly, we get a Cauchy transform for the Rarita-Schwinger operator $R_k:$

\begin{definition}\quad For a domain $\Omega\subset \mathbb{R}^n$ and a function $f: \Omega\times\mathbb{R}^n \longrightarrow Cl_n,$ where $f(x,u)$ is monogenic in $u$, the Cauchy $($or $T_k$-transform$)$ of $f$ is formally defined to be $$(T_kf)(y,v)=-\iint_\Omega \left(E_k(x-y,u,v), f(x,u)\right)_udx^n,\qquad y\in \Omega.$$\end{definition}
\begin{theorem} \quad $R_kT_k\psi=\psi$ for $\psi\in C_0^{\infty}(\mathbb{R}^n,\mathcal{M}_k).$  i.e $$R_k\di\iint_{\mathbb{R}^n}\left(E_k(x-y,u,v), \psi(x,u)\right)_udx^n=\psi(y,v).$$ \end{theorem}
{\bf Proof}\quad For each fixed $y\in \mathbb{R}^n$, let $R(y)$ be a bounded rectangle in $\mathbb{R}^n$ centered at $y.$
Then
$$R_k\di\iint_{\mathbb{R}^n\setminus R(y)}(E_k(x-y,u,v), \psi(x,u))_udx^n=0.$$
Now consider $$\begin{array}{lll}
\di\frac{\partial}{\partial y_i}\di\iint_{R(y)}(E_k(x-y,u,v), \psi(x,u))_udx^n\\
\\
=\di\lim_{\varepsilon\to 0}\di\frac{1}{\varepsilon}\di\iint_{R(y)}(E_k(x-y,u,v)-E_k(x-y-\varepsilon e_i,u,v), \psi(x,u))_udx^n
\end{array}$$
If we translate the rectangle by $\varepsilon$ in $-e_i$ direction, then the derivative will be shifted from $E_k$ to $\psi$. Hence the previous integral becomes
$$\begin{array}{lll}
\di\iint_{R(y)}(E_k(x-y,u,v),\di\frac{ \psi(x,u)-\psi(x+\varepsilon e_i)}{\varepsilon})_udx^n+\\
\\
\di\frac{1}{\varepsilon}\di\iint_{(R(y+\varepsilon e_i)\setminus R(y))\cup (R(y)\setminus R(y+\varepsilon e_i))} (E_k(x-y,u,v), \psi(x,u)-\psi(x+\varepsilon e_i,u))_udx^n\end{array}$$
When $\varepsilon$ tends to zero, the integral is equal to
$$\begin{array}{lll}\di\iint_{R(y)}(E_k(x-y,u,v), \di\frac{\partial \psi(x,u)}{\partial x_i})_udx^n

+\di\int_{\partial R_1(y)\cup \partial R_2(y)}(E_k(x-y,u,v), \psi(x,u))_ud\sigma(x)
\end{array}$$
where $\partial R_1(y)$ and $\partial R_2(y)$ are the two faces of $R(y)$ with normal vectors $\pm e_i.$
So $$\begin{array}{ll}D_y\di\iint_{R(y)}(E_k(x-y,u,v), \psi(x,u))_udx^n\\
\\
=\di\iint_{R(y)}\sum_{i=1}^n e_i(E_k(x-y,u,v), \di\frac{\partial {\psi(x,u)}}{\partial {x_i}})_udx^n\\
\\
+\di\int_{\partial R(y)}n(x)(E_k(x-y,u,v),\psi(x,u))_ud\sigma(x).\end{array}$$
When the volume of $R(y)$ tends to zero, the first integral tends to zero by the homogeneity of the kernel $E_k$. So we shall concentrate attention on the integral
$$P_k\di\int_{\partial R(y)}n(x)(E_k(x-y,u,v),\psi(x,u))_ud\sigma(x).$$
This is equal to
$$P_k\di\int_{\partial R(y)}\di\int_{\mathbb{S}^{n-1}}n(x)E_k(x-y,u,v)\psi(x,u)dS(u)d\sigma(x),$$ which in turn is equal to
$$\begin{array}{ll}
P_k\di\int_{\partial R(y)}\di\int_{\mathbb{S}^{n-1}}n(x)E_k(x-y,u,v)\psi(y,u)dS(u)d\sigma(x)\\
\\
+P_k\di\int_{\partial R(y)}\di\int_{\mathbb{S}^{n-1}}n(x)E_k(x-y,u,v)(\psi(x,u)-\psi(y,u))dS(u)d\sigma(x).\end{array}$$
But the last integral on the right side of the above formula tends to zero as the surface area of $\partial R(y)$ tends to zero. Hence we are left with
$$P_k\di\int_{\partial R(y)}\di\int_{\mathbb{S}^{n-1}}n(x)E_k(x-y,u,v)\psi(y,u)dS(u)d\sigma(x).$$
By Stokes' Theorem this is equal to
$$
P_k\di\int_{\partial B(y,r)}\di\int_{\mathbb{S}^{n-1}}n(x)E_k(x-y,u,v)\psi(y,u)dS(u)d\sigma(x).
$$  In turn this is equal to
$$\begin{array}{ll}
P_k\di\int_{\partial B(y,r)}\di\frac{1}{\omega_nc_k}\di\int_{\mathbb{S}^{n-1}}\di\frac{y-x}{\|x-y\|}\di\frac{x-y}{\|x-y\|^n}Z_k(\di\frac{(x-y)u(x-y)}{\|x-y\|^2},v)\psi(y,u)dS(u)d\sigma(x)\\
\\
=P_k\di\int_{\partial B(y,r)}\di\frac{1}{\omega_n c_k}\di\int_{\mathbb{S}^{n-1}}\di\frac{1}{r^{n-1}}Z_k(\di\frac{(x-y)u(x-y)}{\|x-y\|^2},v)\psi(y,u)dS(u)d\sigma(x).
\end{array}$$
By Lemma 6, the integral becomes
$P_k\di\int_{\mathbb{S}^{n-1}}Z_k(u,v)\psi(y,u)dS(u)=P_k\psi(y,v)=\psi(y,v).\quad \blacksquare$

Now we may establish the intertwining operators for the convolution operator $E_k\star.$ More precisely we shall show that:
 \begin{theorem} If $\psi\in C_0^{\infty}(\mathbb{R}^n, \mathcal{M}_k),$ then $$\begin{array}{ll}J_1(\phi,y)\di\iint_{\mathbb{R}^n}(E_k(x'-y',u,v),\psi(x',v))_vd(x')^n\\
 \\
 =\di\iint_{\mathbb{R}^n}(E_k(x-y,u',w)\tilde J_{-1}(\phi,x),\psi(\phi(x),w))_wdx^n,\end{array}$$ where $x'=\phi(x), y'=\phi(y),  u=\di\frac{(cy+d)u'\widetilde{(cy+d)}}{\|cy+d\|^2}$ and $v=\di\frac{(cx+d)w\widetilde{(cx+d)}}{\|cx+d\|^2}$.

Alternatively, $$J_1(\phi,-)E_k\star \psi=E_k \tilde J_{-1}(\phi,-)\star\psi.$$
\end{theorem}
{\bf Proof}\quad First consider inversion, let $\phi(x)=x^{-1}, \phi(y)=y^{-1}.$ Then $$\begin{array}{lll}\di\iint_{\mathbb{R}^n}(E_k(x^{-1}-y^{-1},u,v),\psi(x^{-1},v))_vd(x^{-1})^n\\
\\
=\di\iint_{\mathbb{R}^n}\int_{\mathbb{S}^{n-1}}E_k(x^{-1}-y^{-1},u,v)\psi(x^{-1},v)dS(v)d(x^{-1})^n\\
\\
=\di\iint_{\mathbb{R}^n}\int_{\mathbb{S}^{n-1}}Z_k\left(u,\di\frac{(x^{-1}-y^{-1})v(x^{-1}-y^{-1})}{\|x^{-1}-y^{-1}\|^2}\right)\di\frac{x^{-1}-y^{-1}}{\|x^{-1}-y^{-1}\|^n}\psi(x^{-1},v)dS(v)d(x^{-1})^n.
\end{array}$$

Since
$$x^{-1}-y^{-1}=-y^{-1}(x-y)x^{-1}=-x^{-1}(x-y)y^{-1}=\di\frac{-x}{\|x\|^2}(x-y)\di\frac{y}{\|y\|^2}=\di\frac{-y}{\|y\|^2}(x-y)\di\frac{x}{\|x\|^2},
$$

$$
\begin{array}{ll}Z_k\left(u,\di\frac{(x^{-1}-y^{-1})v(x^{-1}-y^{-1})}{\|x^{-1}-y^{-1}\|^2}\right)=Z_k\left(u,\di\frac{y(x-y)xvx(x-y)y}{\|x\|^2\|y\|^2\|x-y\|^2}\right)\\
\\
=-\di\frac{y}{\|y\|}Z_k\left(\di\frac{yuy}{\|y\|^2},\di\frac{(x-y)w(x-y)}{\|x-y\|^2}\right)\di\frac{y}{\|y\|},~~\mbox{set}~~ w=\di\frac{xvx}{\|x\|^2}.\end{array}$$
Now the previous integral becomes
$$\begin{array}{ll}
\di\iint_{\mathbb{R}^n}\int_{\mathbb{S}^{n-1}}\di\frac{y}{\|y\|}Z_k\left(\di\frac{yuy}{\|y\|^2},\di\frac{(x-y)w(x-y)}{\|x-y\|^2}\right)\di\frac{y}{\|y\|}y\|y\|^{n-2}G(x-y)x\|x\|^{n-2}\\
\\
\qquad \qquad \qquad\qquad\qquad \qquad \qquad \qquad\qquad\qquad \psi(\phi(x),v)\di\frac{1}{\|x\|^{2n}}dS(v)dx^n\\
\\
=\di\iint_{\mathbb{R}^n}\int_{\mathbb{S}^{n-1}}-y\|y\|^{n-2}Z_k\left(\di\frac{yuy}{\|y\|^2},\di\frac{(x-y)w(x-y)}{\|x-y\|^2}\right)G(x-y)\di\frac{x}{\|x\|^{n+2}}\psi(\phi(x),v)dS(v)dx^n\\
\\
=\di\iint_{\mathbb{R}^n}\int_{\mathbb{S}^{n-1}}-y\|y\|^{n-2}E_k(x-y,u',w))\di\frac{x}{\|x\|^{n+2}}\psi(\phi(x),v)dS(v)dx^n,
\end{array}$$
 where $u'=\di\frac{yuy}{\|y\|^2}, w=\di\frac{xvx}{\|x\|^2}$.
  Then the previous integral is
$$\begin{array}{ll}

\di\iint_{\mathbb{R}^n}\int_{\mathbb{S}^{n-1}}-y\|y\|^{n-2}E_k(x-y,u',w))\di\frac{x}{\|x\|^{n+2}}\psi(\phi(x),\di\frac{xwx}{\|x\|^2})dS(w)dx^n
\end{array}$$
Now multiplying both sides of the equation by $\di\frac{y^{-1}}{\|y\|^{n-2}},$ we obtain
$$\begin{array}{lll}\di\frac{y}{\|y\|^{n}}\di\iint_{\mathbb{R}^n}(E_k(x^{-1}-y^{-1},u,v),\psi(x^{-1},v))_vd(x^{-1})^n\\
\\
=\di\iint_{\mathbb{R}^n}\left(E_k(x-y,u',w)\di\frac{x}{\|x\|^{n+2}},\psi(\phi(x),\di\frac{xwx}{\|x\|^2})\right)_wdx^n,
\end{array}$$ where $u=\di\frac{yu'y}{\|y\|^2}$ and $v=\di\frac{xwx}{\|x\|^2}.$

Next, consider orthogonal transformations. We will apply similar arguments used to establish the equation under inversion. Let $\phi(x)=ax\tilde a$ and $\phi(y)=ay\tilde a,$ where $a\in Pin(n).$ Then $$\begin{array}{lll}\di\iint_{\mathbb{R}^n}(E_k(ax\tilde a-ay\tilde a,u,v),\psi(\phi(x),v))_vd(ax\tilde a)^n\\
\\
=\di\iint_{\mathbb{R}^n}(E_k(a(x-y)\tilde a,u,v),\psi(\phi(x),v))_vd(ax\tilde a)^n\\
\\
=\di\iint_{\mathbb{R}^n}\int_{\mathbb{S}^{n-1}}Z_k(u,\di\frac{a(x-y)\tilde ava(x-y)\tilde a}{\|x-y\|^2})\di\frac{a(x-y)\tilde a}{\|x-y\|^n}\psi(\phi(x),v)dS(v)dx^n\\
\\
=\pm \di\iint_{\mathbb{R}^n}\int_{\mathbb{S}^{n-1}}aZ_k(\tilde aua,\di\frac{(x-y)\tilde ava(x-y)}{\|x-y\|^2})\tilde a\di\frac{a(x-y)\tilde a}{\|x-y\|^n}\psi(\phi(x),v)dS(v)dx^n.
\end{array}$$
Set $w=\tilde ava$, then $v=aw\tilde a$. Hence the integral becomes
$$\begin{array}{ll}
\di\iint_{\mathbb{R}^n}\int_{\mathbb{S}^{n-1}}aZ_k(\tilde aua,\di\frac{(x-y)w(x-y)}{\|x-y\|^2})\di\frac{(x-y)\tilde a}{\|x-y\|^n}\psi(\phi(x),v)dS(v)dx^n\\
\\
=\di\iint_{\mathbb{R}^n}\int_{\mathbb{S}^{n-1}}aE_k(x-y,u',w)\tilde a\psi(\phi(x),aw\tilde a)dS(aw\tilde a)dx^n,
\end{array}$$ where $u'=\tilde aua.$ Now multiplying both sides of the equation by $a^{-1},$ we have
$$\begin{array}{ll}
\tilde a\di\iint_{\mathbb{R}^n}(E_k(ax\tilde a-ay\tilde a,u,v),\psi(\phi(x),v))_vd(ax\tilde a)^n\\
\\
=\di\iint_{\mathbb{R}^n}E_k(x-y,u',w)\tilde a\psi(\phi(x),aw\tilde a)dS(w)dx^n,\end{array}$$
where $u'=\tilde aua$ and $v=aw\tilde a.$ By the Iwasawa decomposition of $\phi(x)=(ax+b)(cx+d)^{-1}$, we obtain
$$\begin{array}{ll}J_1(\phi,y)\di\iint_{\mathbb{R}^n}(E_k(x'-y',u,v),\psi(x',v))_vd(x')^n\\
\\
=\di\iint_{\mathbb{R}^n}(E_k(x-y,u',w)\tilde J_{-1}(\phi,x),\psi(\phi(x),w))_wdx^n,\end{array}$$
where $J_1(\phi,x)=J(\phi,x)=\di\frac{\widetilde{cx+d}}{\|cx+d\|^n}, J_{-1}(\phi,x)=\di\frac{cx+d}{\|cx+d\|^{n+2}}, x'=\phi(x), y'=\phi(y), u=\di\frac{(cy+d)u'\widetilde{(cy+d)}}{\|cy+d\|^2}$ and $v=\di\frac{(cx+d)w\widetilde{(cx+d)}}{\|cx+d\|^2}$. Alternatively, $$J_1(\phi,-)E_k\star\psi =E_k\tilde J_{-1}(\phi,-)\star\psi.\qquad \blacksquare$$


{Charles F. Dunkl}

{Department of Mathematics, University of Virginia, Charlottesville, VA 22904-4137, USA.}
{Email:cfd5z@virginia.edu}

{Junxia Li}

{Department of Mathematics, University of Arkansas, Fayetteville, AR 72701, USA.}
{Email:jxl004@uark.edu }

{John Ryan}

{Department of Mathematics, University of Arkansas, Fayetteville, AR 72701, USA.}
{Email:jryan@uark.edu}

{Peter Van Lancker}

{Faculty of Applied Engineering Sciences, University College of Gent, Member of Gent University, Schoonmeerstaat 52, 9000 Gent, Belgium.}

{Email:Peter.VanLancker@hogent.be}
\end{document}